\documentclass[prx,aps,twocolumn,showpacs]{revtex4-1}
\draft 
\usepackage{graphicx}
\usepackage{amsmath}
\usepackage{amssymb}

\begin{document}

\title{Multiplex congruence network of natural numbers} 

\author{Xiao-Yong Yan$^{1,5}$, Wen-Xu Wang$^{2,6,}$}
\email{wenxuwang@bnu.edu.cn}
\author{Guan-Rong Chen$^{3}$}
\author{Ding-Hua Shi$^{4,}$}
\email[]{shidh2001@shu.edu.cn}
\affiliation{
$^{1}$Systems Science Institute, Beijing Jiaotong University, Beijing 100044, China\\
$^{2}$School of Systems Science, Beijing Normal University, Beijing 100875, China\\
$^{3}$Department of Electronic Engineering, City University of Hong Kong, Hong Kong SAR, China\\
$^{4}$Department of Mathematics, Shanghai University, Shanghai 200444, China\\
$^{5}$Big Data Research Center, University of Electronic Science and Technology of China, Chengdu 611731, China\\
$^{6}$Business School, University of Shanghai for Science and Technology, Shanghai 200093, China\\
}

\date{\today}

\begin{abstract}

Congruence theory has many applications in physical, social, biological
and technological systems. Congruence arithmetic has been a
fundamental tool for data security and computer algebra.
However, much less attention was
devoted to the topological features of congruence relations among
natural numbers. Here, we explore the congruence relations in the
setting of a multiplex network and unveil some unique and
outstanding properties of the multiplex congruence network.
Analytical results show that every layer therein is a sparse and
heterogeneous subnetwork with a scale-free topology.
Counterintuitively, every layer has an
extremely strong controllability in spite of its scale-free
structure that is usually difficult to control. Another amazing
feature is that the controllability is robust against targeted
attacks to critical nodes but vulnerable to random failures, which
also differs from normal scale-free networks. The multi-chain
structure with a small number of chain roots arising from each
layer accounts for the strong controllability and the abnormal feature.
The multiplex congruence network offers a graphical solution to the
simultaneous congruences problem,
which may have implication in cryptography based
on simultaneous congruences. Our work also gains insight into the
design of networks integrating
advantages of both heterogeneous and homogeneous networks
without inheriting their limitations.
 
\end{abstract}

\maketitle 

\section*{Introduction} \label{sec:1}

Congruence is a fundamental concept in number theory. Two integers
$a$ and $r$ are said to be congruent modulo a positive integer $m$
if their difference $a-r$ is integrally divisible by $m$, written
as $a \equiv r \ ({\rm mod} \ m)$ \cite{hardy79}. Congruence
theory has been widely used in physics, biology, chemistry,
computer science, and even music and
business~\cite{n2p92,guterl94,ding96,yan02,schr08}. Because of the
limited computational and storage ability of computers, congruence
arithmetic is particularly useful and applicable to computing
with numbers of infinite length~\cite{yan02}. Significant and
representative applications include generating random
numbers~\cite{knuth97}, designing hash functions~\cite{preneel94}
and checksumming in error detections~\cite{yan02,wang05}. As a
cornerstone of modern cryptography, congruence arithmetic has been
successfully
used in public-key encryption~\cite{salo13}, secret
sharing~\cite{adi79}, digital authentication, and many other data
security applications~\cite{ding96,yan02,schr08}.

\begin{figure*}[ht]
\center
\includegraphics[width=15.0cm]{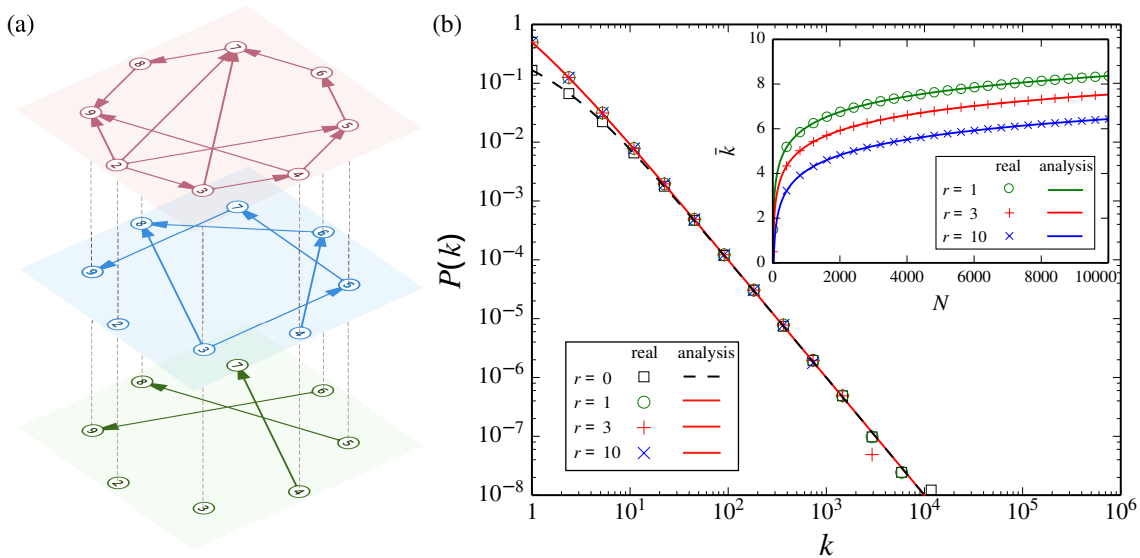}
\caption{{\bf Topology of the multiplex congruence network (MCN).} (a) An MCN with three layers $G(r=\{1,2,3\},9)$, in which each
direct link $L_{ij}$ satisfies the congruence relation $j \equiv r
\ ({\rm mod} \ i)$. (b) Out-degree distributions of congruence
networks with same $N=10000$. The solid lines are the analytical
results from Eq.~(\ref{eq:deg}), and the dashed line is the
analytical result from Eq.~(\ref{eq:deg0}). The insert is the
average degree $\bar{k}$ as a function of $N$ for congruence
networks with different $r$. The solid lines are analytical
results from Eq. (\ref{eq:av}) in the section of \textbf{Methods}. } \label{fig1}
\end{figure*}

Despite the well-established congruence theory with a broad
spectrum of applications,
a comprehensive understanding of the congruence relation among
natural numbers is still
lacking. Our purpose is to uncover some intrinsic properties of
the network consisting of natural numbers with congruence
relations. A link in the congruence network is defined in terms of
the congruence relation
$j \equiv r \ ({\rm mod} \ i)$, where $r$ is the reminder of $j$
divided by $i$. For a fixed value of $r$,
we discern an infinite set of integer pairs $(i,j)$. For each pair of
such integers, a directed link from $i$ to $j$ (suppose
$i<j$) characterizes the congruence relation between $i$ and $j$,
giving rise to a congruence network for a given reminder $r$. Let
$G(r,N)$ denote a
congruence network, where $N$ is the largest integer considered.
Note that congruence networks associated with different values of
$r$ share the same set of nodes (integers), thereby
a multiplex network \cite{yuan14} with a number of layers is formed,
as shown in Fig.~\ref{fig1}(a).

To our knowledge, the multiplex congruence network (MCN) has not been
explored in spite of some effort dedicated to complex networks associated
with natural numbers~\cite{corso04,zhou06,luque08,shi10,garcia14,sheka15}. We will
demonstrate several
unique and prominent properties of the MCN regarding some typical
dynamical processes. Specifically, analytical results will show
that all layers of the MCN are sparse with
the same power-law degree distribution.
A counterintuitive property of the MCN is that every layer of the MCN has an
extremely strong controllability, which significantly differs from
ordinary scale-free networks requiring
a large fraction of driver nodes. To steer the network in a layer,
the minimum number of driver nodes is nothing but the
reminder $r$ that is negligible as compared to the network size.
The controllability of the MCN is also very robust against
targeted removal of nodes but relatively vulnerable to random
failures, which is also in sharp contrast to ordinary scale-free
networks. This amazing
robustness against attacks can be
interpreted in terms of the multi-chain structure in MCN. The MCN
therefore sheds light on the design of heterogeneous networks with
high searching efficiency~\cite{klainberg} and strong
controllability simultaneously.
Another application of the MCN is that it can graphically solve the simultaneous congruences problem in a more
intuitive way than currently used methods, such as the Garner's algorithm. 
The solution of the simultaneous congruences
problem is to locate
common neighbors of relevant numbers in different layers. 
This alternative approach by virtue of the MCN may have implication in
cryptography based on simultaneous congruences.

To our knowledge, the multiplex congruence network (MCN) has not been explored in spite of some effort dedicated to complex networks associated with natural numbers~\cite{corso04,zhou06,luque08,shi10,garcia14,sheka15}. We will
demonstrate several 
unique and prominent properties of the MCN regarding some typical
dynamical processes. Specifically, analytical results will show 
that all layers of the MCN are sparse with 
the same power-law degree distribution.
A counterintuitive property of the MCN is that every layer of the MCN has an 
extremely strong controllability, which significantly differs from
ordinary scale-free networks requiring 
a large fraction of driver nodes. To steer the network in a layer,
the minimum number of driver nodes is nothing but the number of
reminder $r$ that is negligible as compared to the network size.
The controllability of the MCN is also very robust against
targeted removal of nodes but relatively vulnerable to random
failures, which is also in sharp contrast to ordinary scale-free
networks. This amazing 
robustness against attacks can be 
interpreted in terms of the multi-chain structure in MCN. The MCN
therefore sheds light on the design of heterogeneous networks with
high searching efficiency~\cite{klainberg} and strong
controllability simultaneously.
One of the most important applications of the MCN is that it allows 
to solve the simultaneous congruences problem in a much more
efficient way than existing methods using the Gauss algorithm. To
tackle the problem with high computational efficiency, the network
of each layer should be stored in advance, e.g., in a distributed
storage system. The solution of the simultaneous congruences
problem is to locate 
common neighbors of relevant numbers in different layers. This new 
efficient approach by virtue of the MCN offers deeper and broader 
insight into the modern cryptography and has potential
applications in many other 
fields, such as communications and computer science.

\section*{Results} \label{sec:3}

\subsection*{Topology of MCN}

MCN consists of a number of congrence networks (layers)
$G(r>0,N)$, as shown in Fig.~\ref{fig1}(a). Each layer contains
all the natural numbers larger than $r$ but less than or equal to
$N$, so the size (number of nodes) of a layer is $N-r$. The
remainder $r$ is a parameter that determines the structure of
congruence network. When $r=0$, the congruence network reduces to
a divisibility network~\cite{shi10}, in which the dividend links
to all of its divisors except itself.

The out-degrees of nodes are heterogeneous in each layer of the
MCN. We have analytically derived the distribution of the
out-degrees in the thermodynamic limit (see details in the section
of \textbf{Methods}):
\begin{equation}
P(k)= \frac{1}{k(k+1)}.
\end{equation}
For large $k$, the out-degree distribution becomes $P(k) \sim
k^{-2}$, thus $G(r,N)$ is a typical scale-free network. All
$G(r>0,N)$ have similar out-degree distributions, as shown in Fig.
\ref{fig1}(b), but the divisibility network $G(0,N)$ has a
different out-degree distribution. For small $k$, $P(k)$ of
$G(0,N)$ deviates from the other networks $G(r>0,N)$. The main
factor that accounts for the difference lies in that half of the
nodes in $G(0,N)$ have no outgoing links, but in $G(r>0,N)$ there
are only $r$ nodes without outgoing links.

\begin{figure*}[htbp]
\center
\includegraphics[width=17cm]{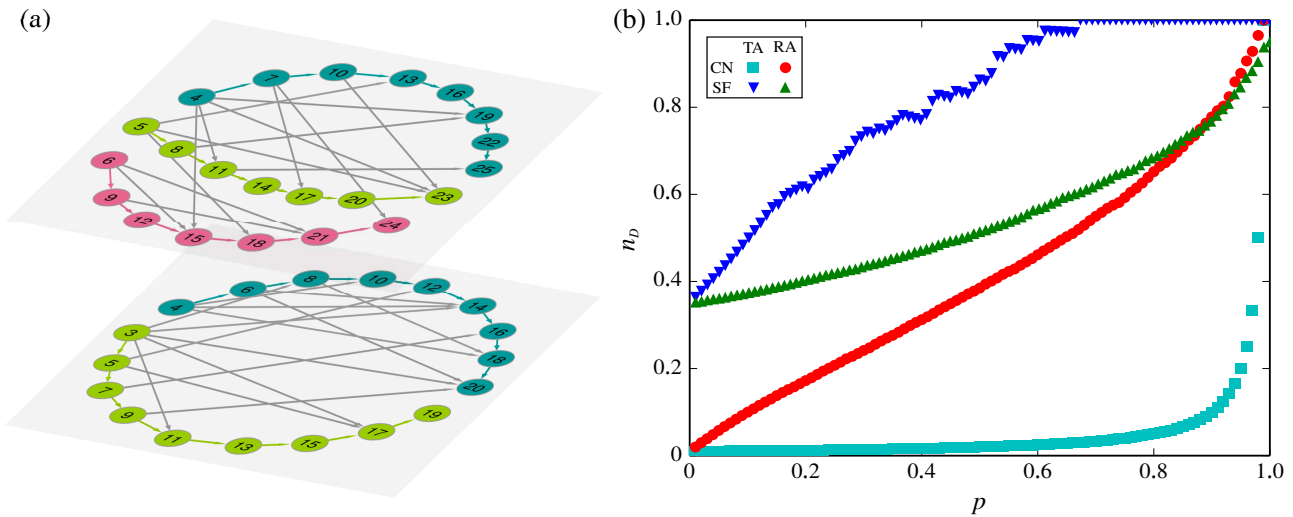}
\caption{ {\bf Controllability of MCN.} (a) Multi-chain structure
in MCN. Different color of nodes and links highlights different
chains in two layers. In each layer, the smaller the node's number,
the larger its out-degree. (b) Driver node density $n_D$ as a function
of the proportion of removed nodes $p$. The nodes are removed
according to two strategies: RA (random attacks: randomly remove
$p$ fraction of nodes) and TA (targeted attacks: remove the top
$p$ fraction of nodes according to their out-degrees). Two
networks are compared: a congruence network (CN) with remainder
$r=1$ and a directed scale-free (SF) network with the same scaling exponent $\gamma
= 2.001$ in both the in- and out-degree distributions. The two networks have the same size $N=100$ and average out-degree $\bar{k}=3.82$. }  \label{fig2}
\end{figure*}

Analytical results demonstrate that the average degree of any
layer increases logarithmically with the network size (see Fig.
\ref{fig1}(b) and the section of \textbf{Methods} for details),
and a larger value of $r$ corresponds to a sparser layer. These
results indicate that $G(r,N)$ is always a sparse network. Hence,
the MCN is compatible with a sparse storage, which is important
for applying the MCN to solve real-world problems.

According to the definition of MCN, the numbers in each layer $G(r > 0,N)$
can be classified into $r$ arithmetic sequences:
\begin{equation}
a_{n}^{i} =i +nr, \label{eq:seq}
\end{equation}
where $1\leq i \leq r, n =1, 2, \cdots, \lfloor \frac{N-i}{r}
\rfloor $ ($\lfloor x \rfloor $ denotes the largest integer not greater than $x$).
The consecutive numbers in the sequence are linked from
small to large, resulting in $r$ chains traversing all nodes in a
layer, as shown in Fig.~\ref{fig2}(a). The root node of a chain is
the minimum number in the chain. There are totally $r$ root nodes
associated with $r$ chains. The end of a chain is always the
maximum number in the chain.
Note that $r=0$ is a special case, because the arithmetic sequence
does not exist in the layer $G(0,N)$, rendering the absence of the
multi-chain structure in the divisibility network.

The above results indicate that although the divisibility network $G(0,N)$ is
a special case of $G(r,N)$, it has some
fundamental differences from $G(r>0,N)$. Only when $r>0$ the
multi-chain structure emerges, and the number of nodes without
outgoing links is negligible. Some evidence has suggested that the multi-chain structure and the absence of nodes with low out-degrees play an
important role in the controllability of complex
networks
\cite{wang12,giulia14}. In the next section, we will further investigate
the controllability properties of the MCN.

\subsection*{Controllability of MCN}

In principle, the MCN composed of natural numbers is not a dynamical system, such that it cannot be controlled. However, because of the multi-chain structure, the MCN provides significant insight into the design of heterogeneous networked systems with strong controllability. Thus, we treat the MCN as a dynamical system and explore its unique and outstanding controllability properties.
The central problem of controlling complex networks is to discern a minimum set of driver nodes, on which external input signals are imposed to fully control the whole system. Let $N_{\rm D}$ denote the minimum number of driver nodes and $n_{\rm D}$ denote the fraction of driver nodes in a network. In general, a network with a smaller value of $n_{\rm D}$ is said to be more controllable.

According to the exact controllability theory for complex networks~\cite{yuan13} and
the sparsity of MCN, we can prove that (see in the section of
\textbf{Methods})
\begin{equation}
N_{\rm D}=r,
\end{equation}
and
\begin{equation}
n_{\rm D}=\frac{r}{N-r}.
\end{equation}
For large $N$ (namely $r \ll N$), $n_{\rm D} \rightarrow 0$ and
$G(r>0,N)$ is considered as highly controllable. Furthermore,
according to both the exact controllability theory~\cite{yuan13}
and the structural controllability theory~\cite{liu11}, the driver
nodes are the $r$ root nodes of the chains in $G(r>0,N)$ (see
details in the section of \textbf{Methods}). Meanwhile, the driver nodes
are the hub nodes with the maximum degree. In comparison, due
to the absence of chains, the divisibility network $G(0,N)$ with
$N_D= \lceil N/2 \rceil$ ($\lceil x \rceil$ denotes the smallest
integer not less than $x$) is hard to control for large $N$.

It has been recognized that scale-free networks are often
difficult to control \cite{liu11}. In particular, Liu et al.
\cite{liu11} have analytically found that when the network size $N
\to \infty$, one must control almost all nodes in order to fully
control a scale-free network with scaling index $\gamma \to 2$.
The MCN is a scale-free network with scaling exponent $\gamma =
2$, but one just needs to control $r$ root nodes to achieve full
control. Such a strong controllability stems from the inherent
multi-chain structure in MCN. From the perspective of structural
controllability, all nodes in the chains in MCN are matched
except the $r$ roots, which need to be controlled. Thus, the MCN is valuable for designing heterogeneous networks with strong controllability.

We also found that the MCN is strongly structurally controllable (SSC) because of the multi-chain structure (see {\bf Methods}), which provides significant insight into the design of heterogeneous and controllable networks without exact link weights. A
network is said to be SSC if and only if its controllability will
not be affected by the link weights in its adjacency matrix~\cite{liu11}, or equivalently, for any distribution of link
weights, the network will be fully controllable from the same set
of driver nodes. 
The SSC property implies that the MCN is robust
against the fluctuation and uncertainty of link weights. This is
an outstanding feature with practical significance
since sometimes link weights are hard to be exactly measured and they are
sometimes time-varying in real situations.

The robustness against attacks is also a significant problem for the design of a controllable networked system~\cite{pu12}. We explore the robustness of the controllability of the MCN against attacks on nodes and find some unique properties, which is useful for the design of practical networks.
On the one hand, due to the
existence of chains rooted in $r$ driver nodes in MCN, targeted
attacks to driver nodes will not destroy the multi-chain structure
nor increase $n_{\rm D}$. Here, nodes critical for
targeted attacks can be identified based on
the rank of node degrees or their hierarchical
structure~\cite{liu12}. In MCN, driver nodes (the $r$ root nodes)
 become such critical nodes. In this regard, MCN is robust
against targeted attacks. On the other hand, random attacks to
nodes may cut some chains. As a result, an additional driver is
required to control each new breakpoint,
leading to an increase of $n_{\rm D}$. Thus, the controllability
of MCN is unusual in resisting
attacks in the sense that it is robust against intentional attacks
but vulnerable to random attacks, which significantly differs from
general scale-free networks. 
The results in comparing the MCN and
SF networks generated by using static model~\cite{goh01} are shown in
Fig.~\ref{fig2}(b). To make an unbiased comparison, a scale-free network with the same scaling exponent $\gamma$ and $\langle k\rangle$ as the MCN is necessary. However, because of the graphicality constraint, it is not possible to generate a random SF network with $\gamma=2$~\cite{del11,baek12}. Thus, we slightly release the requirement by using static model to generate SF networks with the same $\langle k\rangle$ but with $\gamma=2.001$.
Indeed, one can see that $n_{\rm D}$ remains
nearly unchanged under targeted attacks to driver nodes; whereas
random attacks to node causes clear increase of $n_{\rm D}$. This
phenomenon is consistent with our analysis in terms of the
multi-chain structure. Moreover, for the presence of random
attacks, in a wide range of the fractions of failed nodes, the
controllability of MCN is still better than SF networks in
general.

\subsection*{Solving the simultaneous congruences problem}

One of the applications of MCN is that one can
graphically solve the simultaneous congruences problem, which has implication
in communication security and computer
science. In particular, one can find that MCN is exactly a
topological representation of the system of simultaneous
congruences~\cite{schr08}. A system of simultaneous congruences is
a set of congruence equations:
\begin{equation}
\left\{
\begin{array}{c}
x \equiv r_{1} \ ({\rm mod} \ m_{1})\\
x \equiv r_{2} \ ({\rm mod} \ m_{2})\\
\cdots \\
x \equiv r_{k} \ ({\rm mod} \ m_{k})
\end{array}
\right.
\label{eq:sim}
\end{equation}
If the moduli $m_{1},m_{2}\dots m_{k}$ are pairwise coprime, then
a unique solution modulo $m_{1} m_{2} \dots m_{k}$ exists. This is the
\textit{Chinese remainder theorem} (CRT), which has many
applications in computing, coding and cryptography
\cite{ding96,yan02,schr08,sorin07}. A well-known algorithm to
solve the simultaneous congruences in CRT is the Gaussian algorithm
\cite{menezes96}, also known as \textit{Dayanshu} \cite{lib06} in
ancient China. Here, we present an intuitive approach based on MCN
to solve the simultaneous congruences in Eq.~(\ref{eq:sim}).
Firstly, we construct an MCN containing $k$ subnetworks with
remainders $r_{1},r_2 \dots r_k$, respectively.
To focus on the minimum solution of Eq. (\ref{eq:sim}), we set the
maximum number in the network as $N=m_{1} m_{2} \dots m_{k}$.
Then, we find the common successor neighbor of the nodes
$m_{1},m_{2}\dots m_{k}$ in this MCN, which is precisely the
solution $x$.

\begin{figure}[htbp]
\center
\includegraphics[width=9cm]{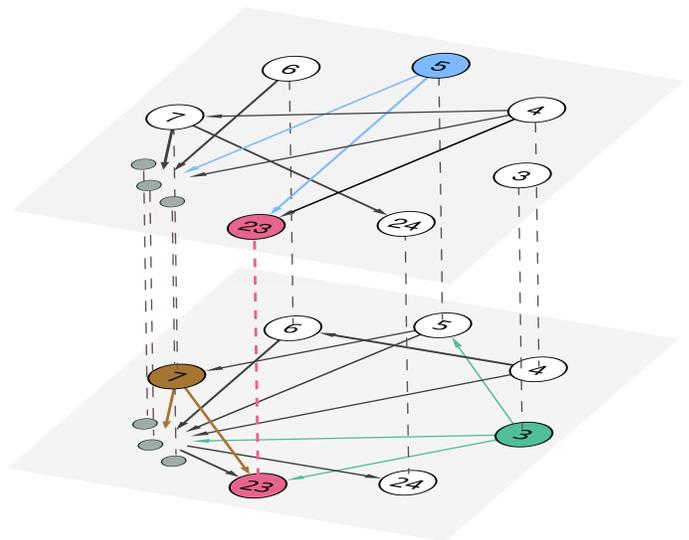}
\caption{ {\bf Solving simultaneous congruences using MCN.} For
visualization, we only show a part of nodes in the MCN. In the
upper layer, the set of successor neighbors of node $5$ is $S_{m=5}=\{8,13,18,23\}$, and similarly $S_{m=3}=\{5,8,11,14,17,20,23\}$
and $S_{m=7}=\{9,16,23\}$ in the lower
layer. Thus the common successor neighbor of the three nodes is
$23$, which is the solution of the simultaneous congruences
problem described by in Eq.~(\ref{eq:sz}).}  \label{fig3}
\end{figure}

We use a well-known example of CRT, recorded in \textit{Sunzi
Suanjing} \cite{lam04}, to demonstrate our approach. The problem
in this example is: `Suppose we have an unknown number of objects.
When grouped in threes, 2 are left out, when grouped in fives, 3
are left out, and when grouped in sevens, 2 are left out. How many
objects are there?' This problem is equivalent to the following
simultaneous congruences
\begin{equation}
\left\{
\begin{array}{c}
x \equiv 2 \ ({\rm mod} \ 3)\\
x \equiv 3 \ ({\rm mod} \ 5)\\
x \equiv 2 \ ({\rm mod} \ 7)
\end{array}
\right.
\label{eq:sz}
\end{equation}
To solve the problem, we first construct an MCN of two layers,
$G(2,105)$ and $G(3,105)$, as shown in Fig. \ref{fig3}. Then, we
find the common successor neighbor of the three moduli, 3, 5 and
7, and finally get the result $x=23$.

It is noteworthy that the traditional algorithm for solving the simultaneous congruences problem, e.g., the Garner's algorithm \cite{menezes96}, is more efficient than our algorithm based on the MCN that is essentially a brute-force search. Thus, it is infeasible to immediately use the graphical approach in data security. However, the graphical algorithm offers new routes to the simultaneous congruences problem in the viewpoint of a complex network, which may be useful to improve the currently used algorithm.

\section*{Discussion} \label{sec:3}
We have defined a multiplex congruence network composed of natural
numbers and uncovered its unique topological features. Analytical
results demonstrate that every layer of the multiplex network is a
sparse and scale-free subnetwork with the same degree
distribution. Counterintuitively, every layer with a scale-free
structure has an extremely strong controllability, which
significantly differs from ordinary scale-free networks. In general,
a scale-free network with power-law degree distribution is harder
to control than homogeneous networks. This is attributed to the
presence of hub nodes, at which dilation arises according to the
structural control theory \cite{liu11}. As a result, downstream
neighbor nodes of hubs are difficult to control. Moreover, due to
a large number of nodes connecting to hubs, scale-free networks
are usually of weak controllability with a large fraction of
driver nodes. In contrast, in spite of the scale-free structure of
the congruence network, the long chains in each layer
considerably inhibit dilation and reduce the number driver nodes.
Furthermore, an interesting finding is that every layer is also
strong structurally controllable in that link weights have no
effect on the controllability. This indicates that the
controllability of the multiplex congruence network is extremely
robust against the inherent limit to precisely accessing link weights
in the real situation. To our knowledge, a scale-free network with
strong structural controllability has not been reported prior to
our congruence network.

An unusual controllability property is that the controllability of
each layer is robust against targeted attacks to driver nodes, but
relatively fragile to random failures of nodes, which is also
different from common scale-free networks. Previously reported
results \cite{pu12,liu12} demonstrate that targeted removal of
high degree nodes and nodes in the top level of a hierarchical
structure causes maximum damage to the network controllability.
Under the two kinds of intentional removals, a network is easier to break
to pieces, such that more driver nodes are required to achieve
full control. Thus, targeted attacks are defined in terms of the two
types of node removals. In the congruence network, high degree
nodes and high level nodes are exactly identical, leading to the
combination of the targeted attacks. Strikingly, the congruence
network is robust against the targeted attacks, because of the
existence of the chains. Targeted attacks will not destroy the
chains in the downstream of the attacked node. As a result, the number of
driver nodes nearly does not increase, even when a large fraction of
nodes has been targeted attacked. The outstanding structural and
controllability properties of the multiplex congruence network are
valuable for designing heterogeneous networks with strong
controllability and high searching efficiency rooted in the scale-free
structure.

Another application of the multiplex congruent network
is to solve the simultaneous congruences problem in a graphical and intuitive manner. 
The multiplex congruence network by converting the algebraic problem of solving simultaneous congruences equations to be a graphical problem of finding common neighbors in a graph, offers an alternative route to the traditional approaches.
Despite this property, the traditional algorithms, such as Gaussian algorithm and Garner's algorithm, outperform the graphical method in computational efficiency. Hence, the graphical method is not applicable in data security at the present. Nevertheless, the graphical approach may inspire the combination of the graphical and algebraic method to improve the current algorithm, which is potentially valuable
in communication security, computer science and many fields relevant
to cryptography. Our work may also stimulate further effort toward
studying of networks arising from natural relationships among
numbers, with outstanding features and applied values. Many
topological insights can be expected from complex networks
consisting of natural numbers.

\section*{Methods}

\subsection*{Deriving the out-degree distribution of MCN}

For a sub-network $G(r>0,N)$ in MCN, the total number of nodes is
$N-r$ and the number of nodes without out-links is $r$. The
out-degree of a node labelled $m$ in the range of
$(\frac{N-r}{2},N-r]$ is 1, because node $m$ can only link to one
node i.e. node $m+r$ in the network; similarly, the out-degree of
a node labelled $m$ in the range of
$(\frac{N-r}{3},\frac{N-r}{2}]$ is 2, because node $m$ can only
link to two nodes i.e. nodes $m+r$ and $2m+r$; similar scenarios
appear to the other nodes. Thus, we can derive the distribution of
out-degrees in the thermodynamic limit, as follows:
\begin{equation}
P(k) =\frac{\frac{N-r}{k}-\frac{N-r}{k+1}}{N-r} = \frac{1}{k(k+1)}, \ k \geq 1.
\label{eq:deg}
\end{equation}

In $G(0,N)$, the numbers larger than $N/2$ have no out-links, i.e.
$P(0)=\frac{N- \lfloor N/2 \rfloor}{N}$, and the numbers in the
range of $(N/3,N/2]$ have only one out-link, i.e.
$P(1)=\frac{\lfloor N/2 \rfloor-\lfloor N/3 \rfloor}{N}$.
Analogously, as $N \to \infty$,
\begin{equation}
P(k) =\frac{\frac{N}{k+1}-\frac{N}{k+2}}{N} = \frac{1}{(k+1)(k+2)}, \ k \geq 0.
\label{eq:deg0}
\end{equation}
This is the same as the in-degree distribution of a growing
network with copying~\cite{ksky05}, which can be regarded as a
random version of the divisibility network~\cite{shi10}.

\subsection*{Calculating the average degree of sparse MCN}

Note that the minimum number $r+1$ in $G(r>0,N)$ has the maximum
out-degree $ \lfloor \frac{N-r} {r+1} \rfloor $ and the second minimum
number $r+2$ has the out-degree $\lfloor \frac{N-r} {r+2} \rfloor $,
and so on. Thus, the average degree $\bar{k}$ of $G(r>0,N)$ is
\begin{equation}
\begin{aligned}
\bar{k} &=\frac{ \sum\limits_{i=r+1}^{N-r}{ \lfloor \frac{N-r}{i} \rfloor } }{N-r}  \\
& \approx \frac{ (N-r) (\ln {(N-r)} +2C-1)-\sum\limits_{i=1}^{r}{ \lfloor \frac{N-r}{i} \rfloor }}{N-r} \\
& \approx \ln(N-r) +H,
\end{aligned}
\label{eq:av}
\end{equation}
where $C$ is the Euler constant and $H = 2C-1
-\frac{\sum\limits_{i=1}^{r}{ \lfloor \frac{N-r}{i} \rfloor
}}{N-r}$ is a constant, which is approximately $C-1-\ln(r)$ when
$r$ is very large.

Eq.~(\ref{eq:av}) is not valid for the divisibility network
$G(0,N)$. In
$G(0,N)$, number 1 has the maximum out-degree $N-1$, and 2 has the
second maximum out-degree $ \lfloor N/2 \rfloor -1$, and so on. In
a similar way to Eq. (\ref{eq:av}), we can obtain the average
degree of $G(0,N)$, as
\begin{equation}
\bar{k} =\frac{\sum\limits_{i=1}^{N}{ \lfloor \frac{N}{i} -1 \rfloor }}{N}  \approx \ln(N) +2C-2,
\label{eq:av0}
\end{equation}
which is consistent with the analytical results of the undirected divisibility network \cite{sheka15}.

\subsection*{Controllability of MCN and identification of driver nodes}
An arbitrary network with linear time-invariant dynamics under control can be described by
\begin{equation}
\dot{\mathbf{x}} = A\mathbf{x} + B\mathbf{u},
\label{eq:ct}
\end{equation}
where the vector $\mathbf{x} = (x_1 , \dots , x_N )^{\mathrm{T}}$
stands for the states of $N$ nodes, $A$ denotes the coupling
matrix (transpose of the adjacency matrix) of a network,
$\mathbf{u} = (u_1 , \dots , u_m )^{\mathrm{T}}$ is the vector of
$m$ input signals, and $B \in \mathbb{R}^{N \times m}$ is the
input matrix.

System (\ref{eq:ct}) is said to be (state) controllable if the
input signal $\mathbf{u}$ imposed on a minimum number $N_{\rm D}$
of driver nodes specified by control matrix $B$ can steer the
state $\mathbf{x}$ from any initial state to any target state in
finite time. The level of controllability of the networked
system~(\ref{eq:ct}) is defined by the fraction $n_{\rm D}$ of
driver nodes in the sense that a complex network is more
controllable if a smaller fraction of driver nodes is needed
to achieve full control. According to Liu {\it et al}.~\cite{liu11}, the key
is to find a matrix $B$ associated with the minimum number of
controllers to ensure full control of system (\ref{eq:ct}).

Because of the sparsity of MCN, $N_{\rm D}$ of a layer $G(r,N)$ is determined by~\cite{yuan13}
\begin{equation}
N_{D} = \max \{1,N-r- \mathrm{rank} (A) \}.
\label{eq:pbh}
\end{equation}
Because in $G(r>0,N)$ each node-labelled number can only link to
the node-labelled numbers larger than itself, the coupling matrix
$A$ is a strictly lower-triangular matrix. Moreover, $A$ of
$G(r>0,N)$ is in a column echelon form because the minimum number
linked from node $m$ is always $m+r$. In other words, the leading
coefficient of the $i$th column in $A$ is precisely in the
$(i+r)$th row. On the other hand, the last $r$ columns of $A$ are
all zeros because the maximum $r$ nodes in $G(r>0,N)$ have no
out-links, so the rank of $A$ is exactly $N-2r$. An example of
matrix $A$ of $G(1,9)$ is
\begin{equation}
A_{G(1,9)}=
\begin{pmatrix}
 0 &   0 &   0 &   0 &   0 &   0 &   0 &   0  \\
  1 &   0 &   0 &   0 &   0 &   0 &   0 &   0  \\
  0 &   1 &   0 &   0 &   0 &   0 &   0 &   0  \\
  1 &   0 &   1 &   0 &   0 &   0 &   0 &   0  \\
  0 &   0 &   0 &   1 &   0 &   0 &   0 &   0  \\
  1 &   1 &   0 &   0 &   1 &   0 &   0 &   0  \\
  0 &   0 &   0 &   0 &   0 &   1 &   0 &   0  \\
  1 &   0 &   1 &   0 &   0 &   0 &   1 &   0  \\
\end{pmatrix}.
\label{eq:mtx1}
\end{equation}
According to Eq. (\ref{eq:pbh}), we can show that the minimum
number of driver nodes $N_{\rm D}$ for $G(r>0,N)$ is $r$ and the value of non-zero elements in $A$ does not affect ${\rm rank}(A)$, which indicate that $G(r>0,N)$ is SSC.

According to the exact controllability theory~\cite{yuan13}, the
control matrix $B$ to ensure full control of the congruence network
$G(r>0,N)$ should satisfy the following condition
\begin{equation}
\mathrm{rank} [-A,B] = N-r.
\label{eq:cond}
\end{equation}
Notice that the rank of the matrix $[-A, B]$ is contributed by the
number of linearly independent rows, hence the input signals
specified via $B$ should be imposed on the linearly dependence
rows in $A$ so as to eliminate all linear correlations
in Eq.~(\ref{eq:cond}). Apparently, the first $r$ rows in the
coupling matrix $A$ of $G(r>0,N)$ are all zero rows (see
Eq.(\ref{eq:mtx1})), hence the $r$ driver nodes that need to be
controlled to maintain full control are just the minimum $r$ nodes
of the congruence network, i.e. the $r$ roots of the chains in the
congruence network (see Fig. \ref{fig2}(a)).

The coupling matrix of the divisibility network $G(0,N)$ is also a
strictly lower-triangular matrix and in a column echelon form, but
the rank of the matrix is $ \lfloor N/2 \rfloor $, because in
$G(0,N)$ the node with labelled number larger than $ \lfloor N/2
\rfloor $ has no out-links, namely, the last $ \lceil N/2 \rceil $
columns of the matrix are all zeros. An example of $G(0,9)$ is
\begin{equation}
A_{G(0,9)}=
\begin{pmatrix}
 0 &   0 &   0 &   0 &   0 &   0 &   0 &   0 &   0  \\
  1 &   0 &   0 &   0 &   0 &   0 &   0 &   0 &   0  \\
  1 &   0 &   0 &   0 &   0 &   0 &   0 &   0 &   0  \\
  1 &   1 &   0 &   0 &   0 &   0 &   0 &   0 &   0  \\
  1 &   0 &   0 &   0 &   0 &   0 &   0 &   0 &   0  \\
  1 &   1 &   1 &   0 &   0 &   0 &   0 &   0 &   0  \\
  1 &   0 &   0 &   0 &   0 &   0 &   0 &   0 &   0  \\
  1 &   1 &   0 &   1 &   0 &   0 &   0 &   0 &   0  \\
  1 &   0 &   1 &   0 &   0 &   0 &   0 &   0 &   0  \\
\end{pmatrix}
.\label{eq:mtx0}
\end{equation}

Therefore, according to Eq.~(\ref{eq:pbh}), we finally obtain the
minimum number of driver nodes of $G(0,N)$ as $N_D=\lceil N/2
\rceil$, indicating that one must control half of the nodes in
order to control the whole divisibility network. Moreover, the value of non-zero elements in $A$ does not affect ${\rm rank}(A)$, which indicate that $G(0,N)$ is SSC. Thus, the MCN composed of $G(0,N)$ and $G(r>0,N)$ is SSC.

\vspace{1em}

{\bf Acknowledgments:} We thank Dr. Zhengzhong Yuan for useful discussions. X.-Y.Y. was
supported by NSFC under Grant Nos. 61304177, 71525002 and the Fundamental
Research Funds of BJTU under Grant No. 2015RC042. W.-X.W. was
supported by NSFC under Grant No. 61573064. 
G.-R.C. was supported by the Hong Kong Research Grants 
Council under the GRF Grant CityU11208515. 
D.-H.S. was supported by NSFC under Grant No. 61174160.


\begin{thebibliography}{10}
\expandafter\ifx\csname url\endcsname\relax
  \def\url#1{\texttt{#1}}\fi
\expandafter\ifx\csname urlprefix\endcsname\relax\def\urlprefix{URL }\fi
\providecommand{\bibinfo}[2]{#2}
\providecommand{\eprint}[2][]{\url{#2}}

\bibitem{hardy79}
\bibinfo{author}{Hardy, G.~H.} \& \bibinfo{author}{Wright, E.~M.}
\newblock \emph{\bibinfo{title}{An Introduction to the Theory of Numbers, 5th
  ed}} (\bibinfo{publisher}{Clarendon Press, Oxford}, \bibinfo{year}{1979}).

\bibitem{n2p92}
\bibinfo{author}{Waldschmidt, M.}, \bibinfo{author}{Moussa, P.},
  \bibinfo{author}{Luck, J.-M.} \& \bibinfo{author}{Itzykson, C.}
\newblock \emph{\bibinfo{title}{From Number Theory to Physics}}
  (\bibinfo{publisher}{Springer, Berlin}, \bibinfo{year}{1992}).

\bibitem{guterl94}
\bibinfo{author}{Guterl, F.}
\newblock \bibinfo{title}{Suddenly, number theory makes sense to industry}.
\newblock \emph{\bibinfo{journal}{Math Horizons}} \bibinfo{pages}{6--8}
  (\bibinfo{year}{1994}).

\bibitem{ding96}
\bibinfo{author}{Ding, C.}
\newblock \emph{\bibinfo{title}{Chinese Remainder Theorem: Applications in
  Computing, Coding, Cryptography}} (\bibinfo{publisher}{World Scientific,
  Singapore}, \bibinfo{year}{1996}).

\bibitem{yan02}
\bibinfo{author}{Yan, S.~Y.}
\newblock \emph{\bibinfo{title}{Number Theory for Computing}}
  (\bibinfo{publisher}{Springer, Berlin}, \bibinfo{year}{2002}).

\bibitem{schr08}
\bibinfo{author}{Schroeder, M.}
\newblock \emph{\bibinfo{title}{Number Theory in Science and Communication:
  With Applications in Cryptography, Physics, Digital Information, Computing,
  and Self-Similarity}} (\bibinfo{publisher}{Springer, Berlin},
  \bibinfo{year}{2008}).

\bibitem{knuth97}
\bibinfo{author}{Knuth, D.~E.}
\newblock \emph{\bibinfo{title}{The Art of Computer Programming}},
  vol.~\bibinfo{volume}{2} (\bibinfo{publisher}{Addison-Wesley, Cambridge},
  \bibinfo{year}{1997}).

\bibitem{preneel94}
\bibinfo{author}{Preneel, B.}
\newblock \bibinfo{title}{Cryptographic hash functions}.
\newblock \emph{\bibinfo{journal}{Eur. T Telecommun.}}
  \textbf{\bibinfo{volume}{5}}, \bibinfo{pages}{431--448}
  (\bibinfo{year}{1994}).

\bibitem{wang05}
\bibinfo{author}{Wang, X.} \& \bibinfo{author}{Yu, H.}
\newblock \bibinfo{title}{How to break md5 and other hash functions}.
\newblock In \emph{\bibinfo{booktitle}{Advances in Cryptology--EUROCRYPT
  2005}}, \bibinfo{pages}{19--35} (\bibinfo{publisher}{Springer, Berlin},
  \bibinfo{year}{2005}).

\bibitem{salo13}
\bibinfo{author}{Salomaa, A.}
\newblock \emph{\bibinfo{title}{Public-key cryptography}}
  (\bibinfo{publisher}{Springer, Berlin}, \bibinfo{year}{2013}).

\bibitem{adi79}
\bibinfo{author}{Shamir, A.}
\newblock \bibinfo{title}{How to share a secret}.
\newblock \emph{\bibinfo{journal}{Communications of the ACM}}
  \textbf{\bibinfo{volume}{22}}, \bibinfo{pages}{612--613}
  (\bibinfo{year}{1979}).

\bibitem{yuan14}
\bibinfo{author}{Yuan, Z.}, \bibinfo{author}{Zhao, C.}, \bibinfo{author}{Wang,
  W.-X.}, \bibinfo{author}{Di, Z.} \& \bibinfo{author}{Lai, Y.-C.}
\newblock \bibinfo{title}{Exact controllability of multiplex networks}.
\newblock \emph{\bibinfo{journal}{New J Phys.}} \textbf{\bibinfo{volume}{16}},
  \bibinfo{pages}{103036} (\bibinfo{year}{2014}).

\bibitem{corso04}
\bibinfo{author}{Corso, G.}
\newblock \bibinfo{title}{Families and clustering in a natural numbers
  network}.
\newblock \emph{\bibinfo{journal}{Phys. Rev. E}} \textbf{\bibinfo{volume}{69}},
  \bibinfo{pages}{036106} (\bibinfo{year}{2004}).

\bibitem{zhou06}
\bibinfo{author}{Zhou, T.}, \bibinfo{author}{Wang, B.-H.},
  \bibinfo{author}{Hui, P.} \& \bibinfo{author}{Chan, K.}
\newblock \bibinfo{title}{Topological properties of integer networks}.
\newblock \emph{\bibinfo{journal}{Physica A}} \textbf{\bibinfo{volume}{367}},
  \bibinfo{pages}{613--618} (\bibinfo{year}{2006}).

\bibitem{luque08}
\bibinfo{author}{Luque, B.}, \bibinfo{author}{Miramontes, O.} \&
  \bibinfo{author}{Lacasa, L.}
\newblock \bibinfo{title}{Number theoretic example of scale-free topology
  inducing self-organized criticality}.
\newblock \emph{\bibinfo{journal}{Phys. Rev. Lett.}}
  \textbf{\bibinfo{volume}{101}}, \bibinfo{pages}{158702}
  (\bibinfo{year}{2008}).

\bibitem{shi10}
\bibinfo{author}{Shi, D.-H.} \& \bibinfo{author}{Zhou, H.-J.}
\newblock \bibinfo{title}{Natural number network and prime number theorem}.
\newblock \emph{\bibinfo{journal}{Complex Systems and Complexity Science}}
  \textbf{\bibinfo{volume}{7}}, \bibinfo{pages}{52--54} (\bibinfo{year}{2010}).

\bibitem{garcia14}
\bibinfo{author}{Garc{\' i}a-P{\' e}rez, G.}, \bibinfo{author}{Serrano, M.~{\'
  A}.} \& \bibinfo{author}{Bogu{\~ n}{\' a}, M.}
\newblock \bibinfo{title}{Complex architecture of primes and natural numbers}.
\newblock \emph{\bibinfo{journal}{Phys. Rev. E}} \textbf{\bibinfo{volume}{90}},
  \bibinfo{pages}{022806} (\bibinfo{year}{2014}).

\bibitem{sheka15}
\bibinfo{author}{Shekatkar, S.~M.}, \bibinfo{author}{Bhagwat, C.} \&
  \bibinfo{author}{Ambika, G.}
\newblock \bibinfo{title}{Divisibility patterns of natural numbers as a complex
  network}.
\newblock \emph{\bibinfo{journal}{Sci. Rep.}} \textbf{\bibinfo{volume}{5}},
  \bibinfo{pages}{14280} (\bibinfo{year}{2015}).

\bibitem{klainberg}
\bibinfo{author}{Kleinberg, J.~M.}
\newblock \bibinfo{title}{Navigation in a small world}.
\newblock \emph{\bibinfo{journal}{Nature}} \textbf{\bibinfo{volume}{406}},
  \bibinfo{pages}{845--845} (\bibinfo{year}{2000}).

\bibitem{wang12}
\bibinfo{author}{Wang, W.-X.}, \bibinfo{author}{Ni, X.}, \bibinfo{author}{Lai,
  Y.-C.} \& \bibinfo{author}{Grebogi, C.}
\newblock \bibinfo{title}{Optimizing controllability of complex networks by
  minimum structural perturbations}.
\newblock \emph{\bibinfo{journal}{Phys. Rev. E}} \textbf{\bibinfo{volume}{85}},
  \bibinfo{pages}{026115} (\bibinfo{year}{2012}).

\bibitem{giulia14}
\bibinfo{author}{Menichetti, G.}, \bibinfo{author}{Dall'Asta, L.} \&
  \bibinfo{author}{Bianconi, G.}
\newblock \bibinfo{title}{Network controllability is determined by the density
  of low in-degree and out-degree nodes}.
\newblock \emph{\bibinfo{journal}{Phys. Rev. Lett.}}
  \textbf{\bibinfo{volume}{113}}, \bibinfo{pages}{078701}
  (\bibinfo{year}{2014}).

\bibitem{yuan13}
\bibinfo{author}{Yuan, Z.}, \bibinfo{author}{Zhao, C.}, \bibinfo{author}{Di,
  Z.}, \bibinfo{author}{Wang, W.-X.} \& \bibinfo{author}{Lai, Y.-C.}
\newblock \bibinfo{title}{Exact controllability of complex networks}.
\newblock \emph{\bibinfo{journal}{Nat. Commun.}} \textbf{\bibinfo{volume}{4}},
  \bibinfo{pages}{2447} (\bibinfo{year}{2013}).

\bibitem{liu11}
\bibinfo{author}{Liu, Y.-Y.}, \bibinfo{author}{Slotine, J.-J.} \&
  \bibinfo{author}{Barab{\'a}si, A.-L.}
\newblock \bibinfo{title}{Controllability of complex networks}.
\newblock \emph{\bibinfo{journal}{Nature}} \textbf{\bibinfo{volume}{473}},
  \bibinfo{pages}{167--173} (\bibinfo{year}{2011}).

\bibitem{pu12}
\bibinfo{author}{Pu, C.-L.}, \bibinfo{author}{Pei, W.-J.} \&
  \bibinfo{author}{Michaelson, A.}
\newblock \bibinfo{title}{Robustness analysis of network controllability}.
\newblock \emph{\bibinfo{journal}{Physica A}} \textbf{\bibinfo{volume}{391}},
  \bibinfo{pages}{4420--4425} (\bibinfo{year}{2012}).

\bibitem{liu12}
\bibinfo{author}{Liu, Y.-Y.}, \bibinfo{author}{Slotine, J.-J.} \&
  \bibinfo{author}{Barab{\'a}si, A.-L.}
\newblock \bibinfo{title}{Control centrality and hierarchical structure in
  complex networks}.
\newblock \emph{\bibinfo{journal}{PLoS ONE}} \textbf{\bibinfo{volume}{7}},
  \bibinfo{pages}{e44459} (\bibinfo{year}{2012}).

\bibitem{goh01}
\bibinfo{author}{Goh, K.-I.}, \bibinfo{author}{Kahng, B.} \&
  \bibinfo{author}{Kim, D.}
\newblock \bibinfo{title}{Universal behavior of load distribution in scale-free
  networks}.
\newblock \emph{\bibinfo{journal}{Phys. Rev. Lett.}}
  \textbf{\bibinfo{volume}{87}}, \bibinfo{pages}{278701}
  (\bibinfo{year}{2001}).

\bibitem{del11}
\bibinfo{author}{Del~Genio, C.~I.}, \bibinfo{author}{Gross, T.} \&
  \bibinfo{author}{Bassler, K.~E.}
\newblock \bibinfo{title}{All scale-free networks are sparse}.
\newblock \emph{\bibinfo{journal}{Phys. Rev. Lett.}}
  \textbf{\bibinfo{volume}{107}}, \bibinfo{pages}{178701}
  (\bibinfo{year}{2011}).

\bibitem{baek12}
\bibinfo{author}{Baek, Y.}, \bibinfo{author}{Kim, D.}, \bibinfo{author}{Ha, M.}
  \& \bibinfo{author}{Jeong, H.}
\newblock \bibinfo{title}{Fundamental structural constraint of random
  scale-free networks}.
\newblock \emph{\bibinfo{journal}{Phys. Rev. Lett.}}
  \textbf{\bibinfo{volume}{109}}, \bibinfo{pages}{118701}
  (\bibinfo{year}{2012}).

\bibitem{sorin07}
\bibinfo{author}{Iftene, S.}
\newblock \bibinfo{title}{General secret sharing based on the chinese remainder
  theorem with applications in e-voting}.
\newblock \emph{\bibinfo{journal}{Electronic Notes in Theoretical Computer
  Science}} \textbf{\bibinfo{volume}{186}}, \bibinfo{pages}{67--84}
  (\bibinfo{year}{2007}).

\bibitem{menezes96}
\bibinfo{author}{Menezes, A.~J.}, \bibinfo{author}{Van~Oorschot, P.~C.} \&
  \bibinfo{author}{Vanstone, S.~A.}
\newblock \emph{\bibinfo{title}{Handbook of Applied Cryptography}}
  (\bibinfo{publisher}{CRC press, Boca Raton}, \bibinfo{year}{1996}).

\bibitem{lib06}
\bibinfo{author}{Ulrich, L.}
\newblock \emph{\bibinfo{title}{Chinese Mathematics in the Thirteenth
  Century}}, vol.~\bibinfo{volume}{1} (\bibinfo{publisher}{Dover Publications,
  New York}, \bibinfo{year}{2006}).

\bibitem{lam04}
\bibinfo{author}{Lam, L.~Y.} \& \bibinfo{author}{Ang, T.~S.}
\newblock \emph{\bibinfo{title}{leeting Footsteps: Tracing the Conception of
  Arithmetic and Algebra in Ancient China}} (\bibinfo{publisher}{World
  Scientific, Singapore}, \bibinfo{year}{2004}).

\bibitem{ksky05}
\bibinfo{author}{Krapivsky, P.~L.} \& \bibinfo{author}{Redner, S.}
\newblock \bibinfo{title}{Network growth by copying}.
\newblock \emph{\bibinfo{journal}{Phys. Rev. E}} \textbf{\bibinfo{volume}{71}},
  \bibinfo{pages}{036118} (\bibinfo{year}{2005}).

\end{thebibliography}
\end{document}